\documentclass{amsart}
\usepackage{amssymb}
\usepackage{eucal}
\usepackage{amsfonts}
\usepackage{epsfig}
\usepackage[all]{xy}
\usepackage[pdftex]{hyperref}
\let\oldcite\cite                                  

\newtheorem{thm}{Theorem}[section]
\newtheorem{cor}[thm]{Corollary}
\newtheorem{lem}[thm]{Lemma}
\newtheorem{prop}[thm]{Proposition}
\theoremstyle{definition}
\newtheorem{defn}[thm]{Definition}
\theoremstyle{remark}
\newtheorem{rem}[thm]{Remark}
\numberwithin{equation}{section} \theoremstyle{remark}
\newtheorem{ex}[thm]{Example}

\newcommand{\X}{\mathcal{X}}
\newcommand{\Y}{\mathcal{Y}}
\newcommand{\PP}{\mathcal{P}}
\newcommand{\OO}{\mathcal{O}}
\newcommand{\LL}{\mathcal{L}}

\newcommand{\G}{\mathcal{G}}
\newcommand{\HH}{\mathcal{H}}

\let\:=\colon

\newcommand{\Ra}{\Rightarrow}
\newcommand{\lra}{\longrightarrow}
\newcommand{\lla}{\longleftarrow}
\newcommand{\llra}[1]{\stackrel{#1}{\lra}}
\newcommand{\llla}[1]{\stackrel{#1}{\lla}}

\newcommand{\mfG}{\mathfrak{G}} \newcommand{\mfH}{\mathfrak{H}}

\newcommand{\mfuG}{\underline{\mathfrak{G}}}

\newcommand{\mcuX}{\underline{\X}}

\newcommand{\cm}{\mathbf{XMod}}
\newcommand{\gp}{\mathbf{Gp}}
\newcommand{\tgp}{\mathbf{2Gp}}
\newcommand{\wtgp}{\mathbf{W2Gp}}
\newcommand{\grst}{\mathbf{grSt}}

\newcommand{\gpd}{\mathbf{Gpd}}
\newcommand{\C}{\mathsf{C}}
\newcommand{\M}{\operatorname{B}}

\newcommand{\Ob}{\operatorname{Ob}}
\newcommand{\Mor}{\operatorname{Mor}}

\newcommand{\Aut}{\underline{\operatorname{Aut}}}
\newcommand{\AUT}{\mathcal{A}ut}
\newcommand{\Hom}{\operatorname{Hom}}
\newcommand{\Pic}{\operatorname{Pic}}
\newcommand{\PGL}{\operatorname{PGL}}
\newcommand{\GL}{\operatorname{GL}}
\newcommand{\Spec}{\operatorname{Spec}}
\newcommand{\id}{\operatorname{id}}
\newcommand{\coker}{\operatorname{coker}}
\newcommand{\im}{\operatorname{im}}

\newcommand{\grPGL}{\mathcal{PGL}}

\def\smashedlongrightarrow{\setbox0=\hbox{$\longrightarrow$}\ht0=1pt\box0}
\def\risom{\buildrel\sim\over{\smashedlongrightarrow}}

\newcommand{\Gm}{\mathbb{G}_m}
\newcommand{\Gmr}[1]{\mathbb{G}_{m,#1}}
\newcommand{\Ao}[1]{\mathbb{A}^{#1}-\{0\}}
 \newcommand{\bbA}{\mathbb{A}}
\newcommand{\bbZ}{\mathbb{Z}} \newcommand{\bbU}{\mathbb{U}}

\begin{document}

\title[Group actions on  algebraic stacks via butterflies]
      {Group actions on algebraic stacks via butterflies}%

\author{Behrang Noohi}%

\begin{abstract} 
  We introduce an explicit method for studying  actions of a group stack $\G$ 
  on an algebraic stack $\X$. As an example, we study in detail the case where 
  $\X=\PP(n_0,\cdots,n_r)$ is a weighted projective stack over an arbitrary base $S$.
  To this end, we give an explicit description of the group stack of automorphisms
  of $\PP(n_0,\cdots,n_r)$, the  {\em weighted projective general linear 2-group}
  $\PGL(n_0,\cdots,n_r)$. As an application, we use a result of Colliot-Th\'el\`ene 
  to show that for  every linear  algebraic group $G$ over an arbitrary base field $k$  
  (assumed to be reductive   if  $\operatorname{char}(k)>0$)   such that 
  $\operatorname{Pic}(G)=0$,   every action  of $G$ on $\PP(n_0,\cdots,n_r)$ lifts 
  to a linear  action of   $G$ on $\bbA^{r+1}$.
\end{abstract}
\maketitle
\section{Introduction}{\label{S:Intro}}


The aim of this work is to propose a concrete method for studying group 
actions on algebraic stacks. Of course, in its full generality this problem could 
already be very difficult in the case of schemes. The case
of stacks has yet an additional layer of difficulty due to the fact that stacks 
have two types of symmetries: 1-symmetries (i.e., self-equivalences) and 
2-symmetries (i.e., 2-morphisms between self-equivalences). 

Studying actions of a group stack $\G$ on a stack $\X$ can be divided into
two subproblems. One, which is of geometric nature,  
is to understand the two types of symmetries alluded to above; these
can be packaged in a group stack $\AUT\X$.  The other, which is of 
homotopy theoretic  nature, is to get a hold of morphisms $\G \to  \AUT\X$.
Here,  a morphism $\G \to \AUT\X$ means a weak
monoidal functor; two morphisms  $f,g \: \G \to  \AUT\X$ that are related by
a monoidal transformation $\varphi \: f \Ra g$ should be regarded as giving 
rise to the ``same'' action.

Therefore, to study actions of $\G$ on  $\X$ one needs to 
understand the group stack $\AUT\X$, the morphisms $\G \to \AUT\X$, 
and also the transformations between such morphisms.
Our proposed method, uses techniques from 2-group theory to 
tackle these problems. It consists of two steps: 

\begin{quote}
1) finding suitable crossed module models for $\AUT\X$ and $\G$; \\
2) using butterflies \cite{Maps, ButterflyI} to give a geometric description of 
morphisms $\G \to  \AUT\X$ and monoidal transformations between them.
\end{quote}

Finding a `suitable'  crossed module model for $\AUT\X$ may not 
always be easy, but we can go about it by choosing a suitable 
`symmetric enough' atlas $X \to \X$. This can be used to find an 
approximation of $\AUT\X$ (Proposition \ref{P:hard}), and if we are 
lucky (e.g., when $\X=\PP(n_0,\cdots,n_r)$) it gives us the whole $\AUT\X$.

Once crossed module models for $\G$ and $\AUT\X$ are found, the butterfly 
method  reduces the action problem to standard problems about group 
homomorphisms and group extensions, which can be tackled using techniques 
from group theory.

\medskip
\noindent{\em Organization of the paper}

\medskip

Sections $\S$\ref{S:2gp}--$\S$\ref{S:Action} 
are devoted to setting up the basic homotopy theory of 2-group actions and 
using butterflies to formulate our strategy for studying actions.  To illustrate 
our method,  in the subsequent sections we apply these ideas to study group 
actions on weighted projective stacks. In  $\S$\ref{S:Recall} we define weighted
projective general linear 2-groups $\PGL(n_0,n_1,\cdots,n_r)$ and prove  (see 
Theorem \ref{T:2-aut}) that they model  $\AUT\PP(n_0,\cdots,n_r)$;  we
prove this over any base scheme $S$, generalizing the case
$S=\Spec\mathbb{C}$ proved in \cite{B-N}:

\begin{thm} Let
   $\AUT\PP_S(n_0,n_1,\cdots,n_r)$  be the group stack
   of   automorphisms of the weighted projective stack
   $\PP_S(n_0,n_1,\cdots,n_r)$ relative to an arbitrary
   base scheme $S$.  Then, there is a natural equivalence of
   group stacks
       $$\grPGL_S(n_0,n_1,\cdots,n_r) \to \AUT\PP_S(n_0,n_1,\cdots,n_r).$$
   Here,  $\grPGL_S(n_0,n_1,\cdots,n_r)$ stands for the group stack associated
   to the crossed module $\PGL_S(n_0,\cdots,n_r)$ (see $\S$\ref{S:Recall} for
   the defnition).
\end{thm} 

We  analyze the structure of $\PGL(n_0,n_1,\cdots,n_r)$ in  detail in 
$\S$\ref{S:Structure}. In Theorem \ref{T:decomposition} we make explicit 
the structure of $\PGL(n_0,\cdots,n_r)=[\mathbb{G}_m \to G]$ by writing 
$G$ as a semidirect product of a reductive part (product of general linear groups) 
and a unipotent part (successive semidirect product of linear affine groups). 

In light of the two step approach discussed above,
Theorems \ref{T:2-aut} and \ref{T:decomposition}  enable us to study actions
of group schemes (or group stacks, for that matter) on weighted
projective stacks in an explicit manner. This is discussed in $\S$\ref{S:Application}. 
We classify actions of a group scheme $G$ on $\PP_S(n_0,n_1,\cdots,n_r)$
in terms of certain central extensions of $G$ by the multiplicative group $\Gm$.
We also describe the stack structure of the corresponding quotient (2-)stacks.
As a consequence, we obtain the following (see Theorem \ref{T:lift}).

\begin{thm}
 Let $k$ be a field and $G$ a connected linear algebraic group over $k$, 
 assumed to be reductive   if  $\operatorname{char}(k)>0$. Let $\X=\PP(n_0,n_1,\cdots,n_r)$
 be a wighted projective stack over $k$. Suppose that  $\operatorname{Pic}(G)=0$.
 Then, every action of $G$ on $\X$ lifts to a linear action of $G$ on $\bbA^{r+1}$.
\end{thm}

In a forthcoming paper, we use the results of this paper (more specifically, 
Theorem \ref{T:2-aut}), together with the results of \cite{Cohomology},  to give 
a complete classification, and explicit construction, of twisted forms of weighted 
projective stacks; these are the weighted analogues of Brauer-Severi varieties.

\tableofcontents

\section{Notation and terminology}{\label{S:NT}}

Our notation for 2-groups and crossed modules is that  of
\cite{Notes} and \cite{Maps}, to which the reader is referred to for
more on  2-group theory  relevant to this work. In particular, we
use mathfrak letters $\mfG$, $\mfH$ for 2-groups or crossed modules.
By a weak 2-group we mean a strict monoidal category $\mfG$ with 
weak inverses (Definition \ref{D:weakdef}). 
If the inverses are also strict, we call $\mfG$ a strict 2-group. 

By a stack we mean a presheaf of groupoids (and not a category
fibered in groupoids) over a Grothendieck site that satisfies the
decent condition. We use mathcal letters $\X$, $\Y$,... for stacks.

Given  a presheaf of groupoids $\X$  over a site, its
stackification is denoted by $\X^a$. We use the same notation for
the sheafification of a presheaf of sets (or groups).

The $m$-dimensional general linear group scheme over $\Spec R$ is 
denoted by $\GL(m,R)$. When $R=\bbZ$, this is abbreviated to $\GL(m)$. 
The corresponding projectivized general linear group scheme is denoted
by $\PGL(m)$; this notation does not conflict with the notation
$\PGL(n_0,n_1,\cdots,n_r)$ for a  weighted projective general linear
2-group ($\S$\ref{S:Recall}) because in the latter case we always
assume $r\geq 1$.

\section{Review of 2-groups and crossed modules}{\label{S:2gp}}

A  {\em strict 2-group} is a group object in the category of
groupoids. Equivalently, a strict 2-group is a strict monoidal groupoid
$\mfG$ in which every object has a strict inverse; that is,
multiplication by an object induces an isomorphism from $\mfG$ onto
itself. A morphism of 2-groups is, by definition,  a strict
monoidal functor.

The weak 2-groups we will encounter in this paper are less weak than the
ones discussed in  \cite{Notes, Maps}. We hope that this change in terminology
is not too confusing for the reader.

\begin{defn}{\label{D:weakdef}}
  A {\em weak  2-group} is a strict monoidal groupoid $\mfG$ in which
  multiplication by an object induces an equivalence of categories
  from $\mfG$ to itself. By a   morphism of weak  2-groups we mean 
  a strict monoidal functor.   By a {\em weak morphism}  we mean a 
  weak monoidal functor. (We will not encounter weak morphisms until 
  later sections.) 
\end{defn}

The set of isomorphism classes of objects in a 2-group
$\mfG$  is denoted by $\pi_0\mfG$; this is a group. The automorphism
group of the identity object $1 \in \Ob\mfG$ is denoted by
$\pi_1\mfG$; this is an abelian group.

Weak 2-groups and strict monoidal functors between them form a
category $\wtgp$ which contains the category $\tgp$ of strict 2-groups 
as a full subcategory.\footnote{Both $\wtgp$ and  $\tgp$
are {\em 2-categories} but we will ignore the 2-morphisms for the
time being and only look at the underlying 1-category.} Morphisms in
$\wtgp$ induce group homomorphisms on $\pi_0$ and $\pi_1$. In other
words, we have  functors $\pi_0, \pi_1\: \wtgp \to \gp$; the functor
$\pi_1$ indeed lands in the full subcategory of abelian groups. A
morphism between weak  2-groups is called an {\em equivalence} if
the induced homomorphisms on $\pi_0$ and $\pi_1$ are isomorphisms.
Note that an equivalence may not have an inverse.

The following lemma is straightforward.

\begin{lem}{\label{L:easy0}}
 Let $f \: \mfH \to \mfG$ be a morphism of weak  2-groups. Then
 $f$, viewed as a morphism of underlying groupoids, is fully
 faithful if and only if $\pi_0f\: \pi_0\mfH \to \pi_0\mfG$ is
 injective and $\pi_1f\: \pi_1\mfH \to \pi_1\mfG$ is an isomorphism.
 It is an equivalence of groupoids if and only if both
 $\pi_0f$ and $\pi_1f$ are isomorphisms.
\end{lem}

A {\em crossed module} $\mathfrak{G}=[\partial \: G_1 \to G_0]$
consists of  a pair of groups $G_0$ and $G_1$, a group homomorphism
$\partial \: G_1 \to G_0$, and a (right) action of $G_0$ on $G_1$,
denoted $-^a$. This action lifts the conjugation action of $G_0$ on the
image of $\partial$ and descends the conjugation action of $G_1$ on
itself. In other words, the following axioms are satisfied:
   \begin{itemize}
   \item  $\forall \beta \in G_1, \forall a\in G_0, \
                 \partial(\beta^a)=a^{-1}\partial(\beta)a$;
    \item  $\forall \alpha,\beta \in G_1, \
                  \beta^{\partial(\alpha)}=\alpha^{-1}\beta\alpha$.              
 \end{itemize}

It is easy to see that the kernel of $\partial$ is a central (in
particular abelian) subgroup of $G_1$; we denote this abelian group
by $\pi_1\mathfrak{G}$. The image of $\partial$ is always a normal
subgroup of $G_0$; we denote the cokernel of $\partial$ by
$\pi_0\mathfrak{G}$. A  morphism  of crossed modules is a pair of
group homomorphisms which commute with the $\partial$ maps and
respect the actions. Such a morphism induces group homomorphisms on
$\pi_0$ and $\pi_1$.

Crossed modules and morphisms between them form a category, which we
denote by $\cm$. We have functors $\pi_0, \pi_1\: \cm \to \gp$; the
functor $\pi_1$ indeed lands in the full subcategory of abelian
groups. A morphism in $\cm$ is said to be an {\em equivalence} if it
induces isomorphisms on $\pi_0$ and $\pi_1$. Note that an
equivalence may not have an inverse.

There is a well-known natural equivalence of categories 
$\tgp \cong \cm$; see \cite{Notes}, $\S$3.3. This equivalence 
respects the functors $\pi_0$ and $\pi_1$. This way, we can think of a
crossed module as a strict 2-group, and vice versa. For this reason, we
may sometimes abuse terminology and  use the term (strict) 2-group for 
an object  which is actually a crossed module; we hope that this will not cause 
any confusion. Note that $\wtgp$ contains  $\tgp$ as a full subcategory.

\section{2-groups over a site and group stacks}

First a few words on terminology. For us a stack is presheaf of
groupoids over a Grothendieck site (and not a category fibered in 
groupoids)   that satisfies the descent condition. This may be
somewhat unusual for algebraic geometers who are used to categories
fibered in groupoids, but it makes the exposition simpler. Of
course, it is standard that this point of view is equivalent to the
one via categories fibered in groupoids. Just to recall how this
equivalence works, to any category fibered in groupoids $\X$ one can
associate a presheaf $\mcuX$ of groupoids over $\C$ which is defined
as follows. By definition, $\mcuX$ is the presheaf that assigns to
an object $U \in \C$ the groupoids
$\mcuX(U):=\Hom(\underline{U},\X)$, where $\underline{U}$ stands for
the presheaf of sets represented by $U$ and $\Hom$ is computed in
the category of stacks over $\C$. Conversely, to any presheaf of
groupoids one associates a category fibered in groupoids defined via
the Grothendieck construction. For more on this we refer the reader to
\cite{Hollander}, especially $\S$5.2.

\subsection{Presheaves of weak  2-groups over a site}
{\label{SS:weakdef}}

Let $\C$ be a Grothendieck site. Let $\wtgp_{\C}$ be the category of
presheaves of weak  2-groups over $\C$; that is, the category of
contravariant functors from $\C$ to $\wtgp$. We define $\tgp_{\C}$
and $\cm_{\C}$ analogously.  There is a natural equivalence of
categories $\tgp_{\C}\simeq \cm_{\C}$. In particular, we can think
of a presheaf of crossed modules as a presheaf  of (strict) 2-groups.
Note that
$\wtgp_{\C}$ contains  $\tgp_{\C}$ as a full subcategory.

Let $\X$ be a presheaf of groupoids over $\C$. To $\X$ we associate
a presheaf of weak  2-groups $\AUT\X \in \wtgp_{\C}$ which
parametrizes auto-equivalences of $\X$. By definition, $\AUT\X$ is
the functor that associates to an object $U$ in $\C$ the weak 
2-group of self-equivalences of $\X_U$, where $\X_U$ is the
restriction of $\X$ to the comma category $\C_U$. (The `comma
category', or the `over category', $\C_U$ is the category of objects
in $\C$ over $U$.) Notice that in the case where $\X$ is a stack,
$\AUT\X$, viewed as a presheaf of groupoids, is also a stack.
Indeed, $\AUT\X$ is almost a group object in the category of stacks
over $\C$. To be more precise, $\AUT\X$ is a  group stack in the sense
of Definition \ref{D:gr} below.

Let $\mfuG \in \wtgp_{\C}$ be a presheaf of weak  2-groups on $\C$.
We define $\pi_0^{pre}\mfuG$ to be the presheaf $U\mapsto
\pi_0\big(\mfuG(U)\big)$, and $\pi_0\mfuG$ to be the sheaf associated to
$\pi_0^{pre}\mfuG$. Similarly, $\pi_1^{pre}\mfuG$ is defined to be
the presheaf $U\mapsto \pi_1\big(\mfuG(U)\big)$, and $\pi_1\mfuG$
to be the sheaf associated to $\pi_1^{pre}\mfuG$.

We define $\pi_0\mfuG$ and $\pi_1\mfuG$ for a presheaf of
crossed modules $\mfuG \in \cm_{\C}$ in a similar manner. The
equivalence of categories between $\tgp_{\C}$ and $\cm_{\C}$
respects $\pi^{pre}_0$, $\pi_0$, $\pi^{pre}_1$ and $\pi_1$. Lemma
\ref{L:easy0} remains valid in this setting if instead of $\pi_0$
and $\pi_1$ we use $\pi^{pre}_0$ and $\pi^{pre}_1$.

\subsection{Group stacks over $\C$}{\label{SS:gr}}

We recall the definition of a group stack from \cite{Breen}. We modify
Breen's definition by assuming that our group stacks are strictly
associative and have strict units. This is all we will need because the 
group stack $\AUT\X$ of self-equivalences of a stack $\X$ (indeed, 
any presheaf of groupoids $\X$) has this property, and that is all we 
are concerned with in this paper.

 \begin{defn}[\oldcite{Breen}, page 19]{\label{D:gr}}
   Let $\C$ be a Grothendieck site. By a  {\em group stack}
   over $\C$ we mean a stack $\G$ that is a strict  monoid object in
   the category of stacks over $\C$ and for which weak inverses
   exist. By a {\em morphism} of group stacks we mean a strict monoidal functor.
   That is, a morphism of stacks that
  strictly respects the monoidal structures. By a {\em weak morphism} we mean a 
  weak monoidal functor.
 \end{defn}

The condition on  existence of weak inverses means that for every $U
\in \Ob\C$ and every object $a$ in the groupoid $\G(U)$,
multiplication by $a$ induces an equivalence of categories from
$\G(U)$ to itself (or equivalently, an equivalence of stacks from
$\X_U$ to itself). This condition is equivalent to saying that, for
every $U \in \Ob\C$, $\X(U)$ is a weak  2-group. More concisely, it
is equivalent to
        $$\G \times \G \llra{(pr,mult)} \G \times \G$$
being an equivalence of stacks.

\begin{rem}
   It is well known that a weak group stack can always be strictified to a strict one.
   So, theoretically speaking, the strictness of monoidal structure in 
   Definition \ref{D:gr} is not restrictive.
   However, given fixed (strict) group stacks $\G$ and $\HH$,  strict morphisms
   $\HH \to \G$ are {\em not} adequate. We will see in the subsequent sections
   that when studying group actions on stacks we can not avoid weak
   morphisms. In this section, however, we will only discuss strict morphisms. 
\end{rem} 

 Let $\grst_{\C}$ be the
category of group stacks and strict morphisms between them (Definition \ref{D:gr});
this is naturally a full subcategory of  $\wtgp_{\C}$.
There are natural functors
   $$\wtgp_{\C} \to \grst_{\C} \ \  \text{and} \ \  \cm_{\C} \to
                              \grst_{\C}.$$
The former is  the stackification functor that sends a
presheaf of groupoids to its associated stack; note that since the
stackification functor preserves products, we can carry over the
monoidal structure from a presheaf of groupoids to its
stackification. The latter functor is obtained from the former by
precomposing with the natural fully faithful functor $\cm_{\C} \to \wtgp_{\C}$
(see the beginning of $\S\ref{SS:weakdef}$). Given a
presheaf of crossed modules $[\partial \: \underline{G_1} \to
\underline{G_0}]$, the associated group stack has as underlying stack
the quotient stack $[\underline{G_0}/\underline{G_1}]$, where
$\underline{G_1}$ acts on $\underline{G_0}$ by multiplication on the
right (via $\partial$).

 \begin{defn}{\label{D:pi}}
   Let $\X$ be a  presheaf of groupoids over $\C$. We define
   $\pi^{pre}\X$ to be the presheaf that sends an object $U$ in
   $\C$ to the set of isomorphism classes in $\X(U)$. We denote the
   sheaf associated to $\pi^{pre}\X$ by $\pi\X$. For   a global
   section $e$ of $\X$, we define $\Aut_{\X}(e)$ to be sheaf
   associated to the presheaf that sends  an object $U$ in $\C$ to
   the group of automorphisms, in the groupoid $\X(U)$, of the object
   $e_U$; note that  when $\X$ is a stack this presheaf is already a
   sheaf and no sheafification is needed.
 \end{defn}
 
Note that when $\G$ is the underlying  presheaf of groupoids
of a presheaf of weak 2-groups $\mfuG \in \wtgp_{\C}$, then
$\pi_0\mfuG= \pi\G$ and  $\pi_1\mfuG=\Aut_{\G}(e)$, where $e$ is 
the identity section of $\G$.

\subsection{Equivalences of group stacks}{\label{SS:equivalence}}

There are two ways of defining the notion of  equivalence between
group stacks. One way is to regard  them as stacks and use the usual
notion of equivalence of stacks. The other is to regard them as
presheaves of weak  2-groups and use $\pi_0$ and $\pi_1$ (see $\S$\ref{SS:weakdef}). 
The next lemma shows that these two definitions agree.

 \begin{lem}{\label{L:isom}}
  Let $\G$ and $\HH$ be group stacks, and let $f \: \HH \to \G$ be a
  morphism of group stacks. Then, the following are equivalent:
    \begin{itemize}
      \item[($i$)] $f$ is an equivalence of stacks;

      \item[($ii$)] The induced maps $\pi_0f \:
        \pi_0^{pre}\HH \to   \pi_0^{pre}\G$ and $\pi_1f \:
        \pi_1^{pre}\HH \to   \pi_1^{pre}\G$ are isomorphisms
        of presheaves of groups;

      \item[($iii$)] The induced maps $\pi_0f \:
        \pi_0\HH \to   \pi_0\G$ and $\pi_1f \:
        \pi_1\HH \to   \pi_1\G$ are isomorphisms
        of sheaves of groups.
    \end{itemize}
 \end{lem}

 \begin{proof}
  The only non-trivial implication is $(iii) \Ra
  (ii)$. In the proof we will use the following standard fact
  from closed model category theory.

  \begin{quote}
     {\bf Theorem} (\oldcite{Hirsch}, Theorem 3.2.13){\bf .} Let
     $\mathcal{M}$ be a closed model category,  $L$ a localizaing
     class of morphisms in $\mathcal{M}$, and $\mathcal{M}_L$  the
     localized model category. Let $\X$ and $\Y$ be fibrant objects
     (i.e., $L$-local objects) in $\mathcal{M}_L$, and let $f \: \Y
     \to \X$ be a morphism in $\mathcal{M}$ that is a weak
     equivalence in the localized model structure $\mathcal{M}_L$
     (that is, $f$ is an $L$-local weak equivalence). Then, $f$
     is a weak equivalence in $\mathcal{M}$.
  \end{quote}

\vspace{0.1in}

  We will apply the above theorem with
  $\mathcal{M}$ being the model structure on the category
  $\gpd_{\C}$ of presheaves of groupoids on $\C$
  in which weak equivalences
  are morphisms that induce isomorphisms (of presheaves of groups) on
  $\pi_0^{pre}$ and $\pi_1^{pre}$, and fibrations are
  objectwise. We take $L$ to be the class of hypercovers.
  The weak  equivalences
  in the localized model structure will then be the ones inducing
  isomorphism (of sheaves of groups) on $\pi_0$ and $\pi_1$.
  The main reference for this is \cite{Hollander}.

  Let us now prove $(iii) \Ra
  (ii)$. It is shown in \cite{Hollander} that
  $\G$ and  $\HH$ are $L$-local objects
  (see $\S$5.2 and $\S$7.3 of [ibid.]).
  By hypothesis, $f$ induces isomorphisms (of sheaves) on $\pi_0$ and
  $\pi_1$, so it is a weak equivalence in the localized
  model structure. Therefore, since $\G$ and  $\HH$ are $L$-local,
  $f$ is already a weak equivalence in the non-localized  model
  structure. This exactly means that $\pi_0f \:
        \pi_0^{pre}\HH \to   \pi_0^{pre}\G$ and $\pi_1f \:
        \pi_1^{pre}\HH \to   \pi_1^{pre}\G$ are isomorphisms
        of presheaves.
 \end{proof}

 \begin{lem}{\label{L:easy3}}
   Let $\X$ be a presheaf of groupoids over $\C$
   and $\varphi \: \X \to \X^a$ its stackification. Then, we
   have the following (see Definition \ref{D:pi} for notation):
  \begin{itemize}
   \item[($i$)] The induced morphism $\pi\X \to \pi(\X^a)$ is
      an isomorphism of sheaves of sets;

   \item[($ii$)] For every global section $e$ of $\X$,
       the natural map $\Aut_{\X}(e) \to \Aut_{\X^a}(e)$ is
       an isomorphism of sheaves of groups.
  \end{itemize}
 \end{lem}

 \begin{proof}
   This is a simple sheaf theory exercise. We include the proof of
   ($i$). Proof of  ($ii$) is similar.

   First we prove that $\pi\varphi \: \pi\X \to \pi(\X^a)$ is
   injective. Let $U \in \Ob\C$, and let $x,y$ be element in
   $\pi\X(U)$ such that $\pi\varphi(x)=\pi\varphi(y)$. We have to
   show that $x=y$. By passing to a cover of $U$, we may assume $x$
   and $y$ lift to objects $\bar{x}$ and $\bar{y}$ in $\X(U)$. We
   will show that there is an open cover  of $U$ over which
   $\bar{x}$ and $\bar{y}$ become isomorphic. Since
   $\varphi(\bar{x})$ and $\varphi(\bar{y})$ become equal in
   $\pi(\X^a)$, there is a cover $\{U_i\}$  of $U$ such that there
   is an isomorphism $\alpha_i \: \varphi(\bar{x}|_{U_i}) \risom
   \varphi(\bar{y}|_{U_i})$ in the groupoid $\X^a(U_i)$, for every
   $i$. By replacing  $\{U_i\}$ with a finer cover, we may assume
   that $\alpha_i$ come from $\X(U_i)$. (More precisely,
   $\alpha_i=\varphi(\beta_i)$, where $\beta_i$ is a morphism in the
   groupoid $\X(U_i)$.) This implies that, for every $i$,
   $\bar{x}|_{U_i}$ and  $\bar{y}|_{U_i}$ are isomorphic as objects
   of the groupoid $\X(U_i)$. This is exactly what we wanted to prove.

   Having proved the injectivity, to prove the surjectivity it is
   enough to show that every object $x$ in $\pi(\X^a)(U)$ is in the
   image of $\varphi$, possibly after replacing $U$ by an open
   cover. By choosing an appropriate cover, we may assume     $x$
   lifts to $\X^a(U)$. Since $\X^a$ is the stackification of $\X$,
   we may assume, after refining our cover,  that $x$ is in the
   image of  $\X(U) \to \X^a(U)$. The claim is now immediate.
 \end{proof}

 \begin{lem}{\label{L:easy4}}
   Let $\mfuG=[\underline{G_1} \to \underline{G_0}]$ be a presheaf
   of crossed modules, and let
   $\G=[\underline{G_0}/\underline{G_1}]$ be the corresponding
   group stack. Then, we have natural isomorphisms of sheaves
   of groups $\pi_i\mfuG \risom \pi_i\G$, $i=1,2$.
 \end{lem}

 \begin{proof}
  Apply Lemma \ref{L:easy3}.
 \end{proof}

\section{Actions of group stacks}{\label{S:Action}}

In this section we  present an interpretation of an action of a group stack
on a stack in terms of butterflies. We begin with the definition of an action.

\begin{defn}{\label{D:action}}
  Let $\X$ be a stack and $\G$ a group stack. By an {\em action} of $\G$ on $\X$ we mean
  a weak morphism $f \: \G \to \AUT\X$. We say two actions $f$ and $f'$ are 
 {\em equivalent} if there is a monoidal transformation $\varphi \: f \to f'$.
\end{defn}

In the case where $G$ is a group (over the base site), 
it is easy to see that our definition of action is equivalent to Definitions 1.3.(i) of \cite{Romagny}. 
A monoidal transformation $\varphi \: f \to f'$ between two such actions
is the same as the structure of a {\em morphism of $G$-groupoids} (in the sense of 
Definitions 1.3.(ii) of \cite{Romagny}) on the identity map $\id_{\X}\: \X \to \X$, where 
the source and the target are endowed with the $G$-groupoid
structures coming from the actions $f$ and $f'$, respectively.

\subsection{Formulation in terms of crossed modules and butterflies}
{\label{SS:classification}}

Butterflies were introduced in \cite{Maps} as a convenient way of  encoding weak morphisms
between 2-groups (rather, crossed modules representing the 2-groups). The theory
was further extended in \cite{ButterflyI} to the relative case (over 
a Grothendieck site). We will use this theory to translate problems about 2-group actions
on stacks to certain group extension problems.

We begin by recalling the definition of a butterfly  (see \cite{Maps}, Definition 8.1
and \cite{ButterflyI}, $\S$4.1.3).

 \begin{defn}{\label{D:butterfly}}
  Let $\mfG=[\varphi \: G_1\to G_0]$ and $\mfH=[\psi \: H_1\to H_0]$
  be crossed modules. By a {\em butterfly} from $\mfH$ to $\mfG$ we
  mean a commutative diagram of groups
   $$\xymatrix@C=8pt@R=6pt@M=6pt{ H_1 \ar[rd]^{\kappa} \ar[dd]_{\psi}
                          & & G_1 \ar[ld]_{\iota} \ar[dd]^{\varphi} \\
                            & E \ar[ld]^{\sigma} \ar[rd]_{\rho}  & \\
                                       H_0 & & G_0       }$$
  in which both diagonal sequences are complexes, and the NE-SW
  sequence, that is, $G_1 \to E \to H_0$, is short exact. We require
  that $\rho$ and $\sigma$ satisfy the following compatibility with
  actions. For every $x \in E$, $\alpha \in G_1$, and $\beta \in
  H_1$,
     $$\iota(\alpha^{\rho(x)})=x^{-1}\iota(\alpha) x, \ \
          \kappa(\beta^{\sigma(x)})=x^{-1}\kappa(\beta) x.$$
  A {\em morphism} between two butterflies $(E,\rho,\sigma,\iota,\kappa)$
  and $(E',\rho',\sigma',\iota',\kappa')$ is a morphism $f \: E \to E'$
  commuting with all four maps (it is easy to see that such an $f$ is
  necessarily an isomorphism). We define $\M(\mfH,\mfG)$ to be
  the groupoid of butterflies from $\mfH$ to $\mfG$.
 \end{defn}

This definition is justified by the following result
(see \cite{Maps} and \cite{ButterflyI}).

\begin{thm}{\label{T:weakmorphism}}
  Let $\mfG=[G_1\to G_0]$ and $\mfH=[H_1\to H_0]$ be crossed modules
  of sheaves of groups over the cite $\C$. Let $\G=[G_0/G_1]$ and
  $\HH=[H_0/H_1]$ be the corresponding quotient group stacks. Then,
  there is a natural equivalence of groupoids
    $$\M(\mfH,\mfG)\cong \Hom_{\operatorname{weak}}(\HH,\G).$$
  Here, the right hand side stands for the groupoid whose objects are weak
  morphisms of group stacks and whose morphisms are monoidal transformations.
\end{thm}

The above result can be interpreted as follows. A butterfly as in the theorem
gives rise to a canonical zigzag in $\cm_{\C}$
  $$\mfH \llla{\sim} \mathfrak{E} \to \mfG,$$
where $\mathfrak{E} = [H_1\times G_1 \llra{\kappa\cdot\iota} E]$. 
After passing to the associated stacks, it gives rise to a zigzag in $\grst_{\C}$
   $$\HH \llla{\sim} \mathcal{E} \to \G,$$
which after inverting the left map (as a weak morphism), results in a weak morphism 
$\HH \to \G$. It follows from this description of a butterfly that $\pi_0$ and $\pi_1$ 
are functorial with respect to butterflies. Furthermore, the equivalence of Theorem 
\ref{T:weakmorphism} respects $\pi_0$ and $\pi_1$ (see Lemma \ref{L:easy4}).

When $[G_1\to G_0]$ is a crossed module model for the
group stack $\AUT\X$ of auto-equivalences of a stack $\X$, then it follows from 
Theorem \ref{T:weakmorphism} that an action of $\HH=[H_0/H_1]$ on $\X$ is the 
same thing as an isomorphism
class of a butterfly as in Definition \ref{D:butterfly}. In other
words, to give an action of $\HH$ on $\X$, we need to find an
extension $E$ of $H_0$ by $G_1$, together with group homomorphisms
$\kappa \: H_1 \to E$ and $\rho \: E \to G_0$ satisfying the
conditions of Definition \ref{D:butterfly}. This summarizes our strategy  
for studying group actions on stacks. To show its usefulness,
in the subsequent sections we will apply
this method to the case where $\X=\PP_S(n_0,n_1,\cdots,n_r)$
is a weighted projective stack over a base scheme $S$.

\section{Weighted projective general linear 2-groups  $PGL(n_0,n_1,...,n_r)$}{\label{S:Recall}}

In this section we introduce weighted projective general linear 2-group schemes and prove that
they model self-equivalences of weighted projective stacks (Theorem \ref{T:2-aut}).

We begin by some general observations about automorphism 2-groups of  quotient stacks.
From now on, we assume that $\C=\mathbf{Sch}_S$ is the
big site of schemes over a base scheme $S$, endowed with a
subcanonical topology (say, \'etale, Zariski, fppf, fpqc, etc.). 

\subsection{Automorphism 2-group of a quotient stack}{\label{SS:quotient}}

We define a {\em crossed module in $S$-schemes} $[\partial \: G_1 \to
G_0]$ to be  a pair of $S$-group schemes $G_0$ and $G_1$, an $S$-group
scheme homomorphism $\partial \: G_1 \to G_0$, and a (right) action
of $G_0$ on $G_1$ satisfying the axioms of a crossed module. These
are precisely the representable objects in $\cm_{\mathbf{Sch}_S}$;
in other words, a crossed module in schemes $[\partial \: G_1 \to
G_0]$ gives rise to a presheaf of crossed modules
    $$U \mapsto [\partial(U) \: G_1(U) \to G_0(U)].$$

We often abuse terminology and call a crossed module in schemes
over $S$ simply  a {\em strict 2-group scheme over $S$}.

The following two propositions generalize Lemma  8.2 of  \cite{B-N}.

 \begin{prop}{\label{P:hardold}}
   Let $S$ be a base scheme. Let $A$ be an abelian affine group scheme over
   $S$ acting on a $S$-scheme $X$, and let $\X=[X/A]$ be the
   quotient stack. Let $G$ be  those automorphisms of $X$ which
   commute with the $A$ action; this is a sheaf of groups
   on $\mathbf{Sch}_S$. We have the following:

   \begin{itemize}
     \item[(i)] With the trivial action of $G$ on $A$,
       the natural map $\varphi \: A \to G$ becomes a
       crossed modules in $\mathbf{Sch}_S$-schemes.

     \item[(ii)] Let $\G$ be the group stack
       associated to $[\varphi\: A \to G]$. Then, there is a natural
       morphism of group stacks $\G \to \AUT\X$. Furthermore,
       this morphism induces an isomorphism of sheaves of groups
       $\pi_1\G \risom \pi_1(\AUT\X)$.
   \end{itemize}
 \end{prop}

 \begin{proof}
     Part (i) is straightforward, because $\varphi$ maps $A$
     to the center of $G$.  Let $\mfuG$ denote the presheaf of 2-groups associated to 
      $[\varphi\: A \to  G]$.   To prove part (ii), it is enough to 
      construct a morphism of presheaves of
      2-groups $\mfuG \to \AUT\X$ and show that it has the required
      properties.  Stackification of this map gives us the desired map
      (Lemma \ref{L:easy4}).
 
      Let us construct the morphism  $\mfuG \to \AUT\X$. We give the effect of this 
      morphism on the sections over $S$. Since everything commutes
      with base change, the same construction works for every $U \to S$ in the site
      $\mathbf{Sch}_S$ and gives rise to the desired morphism.
      of presheaves.
      
      To define $\mfuG(S) \to \AUT\X(S)$, recall the explicit description
      of the $S$-points of the quotient stack $[X/A]$:

      {\small
      $$\Ob[X/A](S)=\left\{
      \begin{array}{rcl} (T,\alpha) &\vert& T \  \text{an $A$-torsor over $S$} \\
                             & & \alpha \: T \to X  \ \text{an $A$-map}
                                              \end{array}\right\}$$
      $$\Mor[X/H](S)((T,\alpha),(T',\alpha'))=\{f \: T \to T'
        \ \text{an $A$-torsor map s.t.} \ \alpha'\circ f=\alpha\}$$}

  Any element of $g \in G(S)$ induces
  an automorphism of $\X$ relative to $S$ (keep the same
  torsor $T$ and compose $\alpha$ with the action of $g$ on $X$).
  Also, for any element $a \in A(S)$, there is a natural
  2-isomorphism  from the identity automorphism of $\X$ to the automorphism
  induced by  $\varphi(a) \in G(S)$ (which is by
  definition the same as the action of $a$). It is given by the multiplication
  action of $a^{-1}$ on the torsor $T$ (remember that $A$ is abelian)
  which makes the following triangle commute
       $$\xymatrix{ T \ar[r]^{\alpha} \ar[d]_{a^{-1}} &  X  \\
                    T \ar[ru]_{a\circ\alpha}    &           }$$

  Interpreted in the language of 2-groups, this gives a
   morphism of 2-groups   $\mfuG(S) \to \AUT\X(S)$.

  To  prove that  $\G \to \AUT\X$ induces an isomorphism on $\pi_1$, we show that,
  for every   $U \to S$ in the site $\mathbf{Sch}_S$, the morphism of 2-groups  
  $\mfuG(U) \to \AUT\X(U)$ induces an isomorphism on $\pi_1$. Again, we may 
  assume that $U=S$.   We know that the group of 2-isomorphisms
  from the identity automorphism of $\X$ to itself  is naturally isomorphic
  to the group of global sections of the  inertia stack of $\X$. In the
  case $\X=[X/A]$, this is naturally isomorphic to the  group of elements
  of $A(S)$ which act trivially on $X$.   Note that this group is naturally
  isomorphic to $\pi_1\mfuG(S)$. Therefore, the map
  $\mfuG(S) \to \AUT\X(S)$ induces an isomorphism on $\pi_1$. 
\end{proof}

   \begin{prop}{\label{P:hard}}
      Notation being as in Proposition \ref{P:hardold}, assume that $\X$ is a 
      proper Deligne-Mumford stack   over $S$, and that $X \to S$  has geometrically
       connected and reduced fibers. Also, assume that $A$ fits in an extension
             $$0 \to A_0 \to A \to A/A_0\to 0$$
        where $A/A_0$ is finite over $S$ and  $A_0$ has geometrically connected fibers 
        (this is automatics, for example, in the case where $A$ is smooth and the number 
        of  its geometric connected components is a locally constant function on $S$).   
        Then, $\G \to \AUT\X$ is  fully faithful (as a morphism of presheaves of groupoids).
       In particular, the induced map  $\pi_0\G \to \pi_0(\AUT\X)$ of sheaves of groups is
       injective.
   \end{prop}

  \begin{proof}
  As in Proposition \ref{P:hardold}, let $\mfuG$ denote the presheaf of 
  2-groups associated to  $[\varphi\: A \to  G]$.
  We need to show that, for every $U \to S$ in the site $\mathbf{Sch}_S$,
  $\mfuG(U) \to \AUT\X(U)$ is fully faithful; since  $\AUT\X$ is a stack,
  it would then follow that the stackified morphism $\G \to \AUT\X$ is also
  fully faithful.
  
  We may assume that $U=S$. By Proposition \ref{P:hardold}.(ii) and Lemma 
  \ref{L:easy0},   it is enough to prove that  if the action of  $g \in G(S)$ on 
  $\X$ is 2-isomorphic to the identity, then $g$ is of the form 
  $\varphi(a)$, for some $a \in A(S)$. Let us fix such a 2-isomorphism. 
  The effect of this 2-isomorphism  on the $A$-torsor on $X$  corresponding 
  the point $X \to [X/A]$, viewed as an object  in the  groupoid $[X/A](X)$ of 
  $X$-points of  $[X/A]$, is given by   an $A$-torsor map 
  $F \: A\times_S X \to A\times_S X$ which makes
  the following $A$-equivariant triangle commute:
  $$\xymatrix{ A\times_S X \ar[r]^(0.62){\mu} &  X  \\
                    A\times_S X \ar[ru]_(0.55){g\circ\mu}  \ar[u]^{F}    &       }$$
  Here, $A\times_S X$ is the trivial $A$-torsor on $X$
  and $\mu$ is the action of $A$ on $X$.

   Precomposing $F$ with the canonical section $X \to A\times_S X$ (corresponding 
   to the identity element of $A$) and then  projecting onto the first factor, we obtain
   a map $f \: X \to A$ relative to $S$. The proposition follows from the following.

\medskip
   
\noindent{\em Claim.}   The map $f$ is constant, in the sense that it factors through an
   $S$-point $a\: S \to A$ of $A$. Furthermore, the effect of $a$ on $X$ (induced from 
   the action of   $A$ on $X$) is the same as the effect of $g$ on $X$.

\medskip
   Let us prove the claim.  It follows from the commutativity of the
   above diagram that, for any point $x$ in $X$, the effect of $g$ on $x$
   is the same as the effect of $f(x)$ on $x$.\footnote{When we say a ``point'' of $X$
   we mean a scheme $T$ over $S$ and a morphism $T \to X$ relative to $S$.}
   In other words, $f(x)g^{-1}$   leaves $x$ fixed. Applying this to $ax$ instead of $x$, 
   and using the fact that   $a$ and $f(x)g^{-1}$ commute, we find that $f(ax)g^{-1}$ 
   also leaves $x$   fixed, for every $a \in A$.   This implies that, for any point $x$ 
   of $X$, and any $a \in A$,    the element   $r(a,x):=f(ax)f(x)^{-1}$ leaves $x$ fixed. 
   Therefore,  the map     $\rho \:A\times_S X \to A\times_S X$,   $\rho(a,x):=(r(a,x),x)$  
   factors through the stabilizer   group scheme $\tau\: \Sigma \to X$. Thus, we have a 
   commutative triangle
       $$\xymatrix{ A\times_S X \ar[rr]^(0.62){\rho} \ar[rd]_{pr_2} &&  \Sigma  \ar[ld]^{\tau}\\
                  &   X    &           }$$
    
   Now, consider the short exact sequence
         $$0 \to A_0 \to A \to A/A_0\to 0,$$
   where $A_0$ is a group scheme over $S$ with geometrically connected fibers and 
   $A/A_0$ is finite over $S$.   Since $\tau \: \Sigma \to X$    has discrete fibers  
   (because $\X$ is Deligne-Mumford)   the restriction of $\rho$ to $A_0 \times_S X$ 
   factors through the identity section.   Hence, for every $a \in A_0$ and 
   $x \in X$ (over the same point in $S$), $r(a,x)=f(ax)f(x)^{-1}$ is the  identity element 
   of $A$. This implies that $f \: X \to A$ is $A_0$-equivariant (for the trivial action of
   $A_0$ on $A$). So, we obtain   an induced map $\lambda \: [X/A_0] \to A$ (relative 
   to $S$).    Since  $[X/A_0]$ is   finite over $[X/A]$, and $[X/A]$ is proper over $S$, 
   the structure  map $\pi \: [X/A_0] \to S$ is proper. From our assumptions we have that
   $\pi$  has geometrically connected and reduced fibers. Base change then implies 
   that $\pi_*{\mathcal{O}_{[X/A_0]}}=\mathcal{O}_S$. Since $A$ is affine over 
   $S$, it follows that $\lambda$ is constant, i.e.,  factors through a section $a \: S \to A$. 
   Since   $f \: X \to A$ factors through $\lambda$, it also factors through $a$. By 
   construction, the effect of  $a$ on $X$ is the same as  the effect of $g$ on $X$, which
   is what we wanted to prove.
 \end{proof}

\subsection{Weighted projective general linear 2-groups}{\label{SS:wpgl}}
Since the construction of the weighted projective stacks, and also
of the weighted projective general linear 2-group schemes,
commutes with base change, we can  work over
$\mathbb{Z}$. We begin with some notation. We denote the multiplicative group scheme
over $\Spec \bbZ$ by $\Gmr{\bbZ}$, or simply $\Gm$. The affine
($r+1$)-space over a base scheme $S$ is denoted by $\bbA_S^{r+1}$;
when the base scheme is $\Spec R$ it is denoted by $\bbA_R^{r+1}$,
and when the base scheme is $\Spec \bbZ$ simply by
$\mathbb{A}^{r+1}$. Since $r$ will be fixed throughout this
section, we will usually denote $\bbA^{r+1}_S-\{0\}$ by $\bbU_S$.
We will abbreviate $\bbU_{\Spec R}$ and $\bbU_{\Spec \bbZ}$ to $\bbU_R$ 
and $\bbU$, respectively. We fix a Grothendieck topology on $\mathbf{Sch}_S$ that is not
 coarser than Zariski.

 Let $n_0,n_1,\cdots,n_r$ be a sequence of positive integers, and
 consider the weight $(n_0,n_1,\cdots,n_r)$ action of $\Gm$ on
 $\bbU=\Ao{r+1}$. (That is, for every  scheme $T$, an element $t \in \Gm(T)$ acts
 on $\bbU_T$ by multiplication by $(t^{n_0},t^{n_1},\cdots,t^{n_r})$.) 
 The quotient stack of this action is called the {\em weighted projective stack} of weight
 $(n_0,n_1,\cdots,n_r)$ and is denoted by $\PP_{\bbZ}(n_0,n_1,\cdots,n_r)$, 
 or simply by $\PP(n_0,n_1,\cdots,n_r)$. The {\bf weighted projective general 
 linear 2-group scheme} $\PGL(n_0,n_1,\cdots,n_r)$ is defined
 to be the 2-group scheme associated to the crossed module
    $$[\partial \: \Gm \to G_{n_0,n_1,\cdots,n_r}],$$
 where $G_{n_0,n_1,\cdots,n_r}$ is the group scheme, over $\bbZ$,
 of all $\Gm$-equivariant (for the above weighted
 action) automorphisms of $\bbU$. More precisely, the
 $T$-points of $G_{n_0,n_1,\cdots,n_r}$ are automorphisms
    $$f \: \bbU_T  \to  \bbU_T $$
 that commute with the $\Gm$-action.  The homomorphism 
 $\partial \: \Gm \to G_{n_0,n_1,\cdots,n_r}$ is the one induced from the
 $\Gm$-action itself. We take the action of $G_{n_0,n_1,\cdots,n_r}$
 on $\Gm$ to be  trivial. The associated group stack is denoted by
 $\grPGL(n_0,n_1,\cdots,n_r)$, and is called the {\em projective general linear} 
 group stack of weight $(n_0,n_1,\cdots,n_r)$.

 The following theorem says that a weighted projective general
 linear 2-group scheme is a  model for the group stack of self-equivalences
 of the corresponding weighted projective stack. A special case of
 this theorem (namely, the case where the base scheme is  $\mathbb{C}$)
  was proved in (\cite{B-N}, Theorem 8.1).   We briefly sketch how the proof in [ibid.] 
  can be modified to cover the general case.

 \begin{thm}{\label{T:2-aut}} Let
   $\AUT\PP(n_0,n_1,\cdots,n_r)$  be the group stack
   of   automorphisms of the weighted projective stack
   $\PP(n_0,n_1,\cdots,n_r)$.  Then, the natural map
       $$ \grPGL(n_0,n_1,\cdots,n_r) \to \AUT\PP(n_0,n_1,\cdots,n_r)$$
   is an equivalence of group stacks. In particular, we have
   isomorphisms of sheaves of groups
  {\small
       $$\hspace{-0.32in} \pi_0\AUT\PP(n_0,n_1,\cdots,n_r) \cong
          \pi_0\grPGL(n_0,n_1,\cdots,n_r)\cong \pi_0\PGL(n_0,n_1,\cdots,n_r),$$
       $$ \pi_1\AUT\PP(n_0,n_1,\cdots,n_r) \cong
          \pi_1\grPGL(n_0,n_1,\cdots,n_r) \cong \pi_1\PGL(n_0,n_1,\cdots,n_r) \cong\mu_d,$$
  }
    where $d=\gcd(n_0,n_1,\cdots,n_r)$ and $\mu_d$ stands for the multiplicative group
    scheme of $d^{th}$ roots of unity.
 \end{thm}

In order to prove our main result (Theorem \ref{T:2-aut}) we need the following result
about line bundles on weighted projective stacks. For more details on this, the reader 
is referred to \cite{Picard}. More general results about Picard stacks of algebraic
stacks can be found in \cite{Brochard}.

\begin{prop}{\label{P:linebundles}}
 Let $\PP=\PP_S(n_0,n_1,\cdots,n_r)$, where $S=\Spec R$
 is the spectrum of a local ring. Then every line bundle on
 $\PP$ is of the form $\OO(d)$ for some $d\in \bbZ$.
\end{prop}

\begin{proof}
  In the proof we use stack versions of Grothendieck's base change and
  semicontinuity results (\cite{Hart}, III. Theorem 12.11). We will assume 
  that $R$ is Noetherian. 

In the case where $R$ is a field, the assertion is easy to prove
using the fact that the Picard group of $\PP$ is isomorphic to the
Weil divisor class group. To prove the general case, let $x$ be
the closed point of $S=\Spec R$. Let $\LL$ be a line bundle on
$\PP$. After twisting with on appropriate $\OO(d)$, we may assume
$\LL_x\cong\OO$.   We will show that $\LL$ is trivial. We have
$H^1(\PP_x,\LL_x)=H^1(\PP_x,\OO_x)=0$. Hence, by semicontinuity,
$H^1(\PP_y,\LL_y)=0$ for every point $y$ of $S$. Base change
implies that $R^1f_*(\LL)=0$, and that $R^0f_*(\LL)=f_*(\LL)$ is
locally free (necessarily of rank 1). Therefore, $f_*(\LL)$ is
free of rank 1 and, by base change, $H^0(\PP_y,\LL_y)$ is
1-dimensional as a $k(y)$-vector space, for every $y$ in $S$. In
fact, this is true for every tensor power $\LL^{\otimes n}$, $n
\in \bbZ$. So, $\LL_y$ is trivial for every $y$ in $S$. (Note
that, when $k$ is a field, $\dim_kH^0(\PP_k(n_0,n_1,\cdots,n_r),\OO(d))$ 
is equal to the number of solutions of the equation 
$a_1n_0+a_2n_1+\cdots+a_rn_r=d$ in non-negative integers $a_i$.)

Now let $s$ be a generating section of $f_*(\LL)\cong R$. It
follows that $f^*(s)$ is a generating section of $\LL$. So $\LL$
is trivial.
\end{proof}

\begin{proof}[Proof of Theorem \ref{T:2-aut}]
   We apply Propositions \ref{P:hardold} and  \ref{P:hard} with $S=\Spec{\bbZ}$,
   $X=\Ao{r+1}$, and $H=\Gm$. This  implies that
   $$ \grPGL(n_0,n_1,\cdots,n_r) \to \AUT\PP(n_0,n_1,\cdots,n_r)$$
   is a fully faithful morphism of stacks. That is, for every
   scheme $U$, the morphism of groupoids 
   $$\grPGL(n_0,n_1,\cdots,n_r)(U) \to
   \AUT\PP(n_0,n_1,\cdots,n_r)(U)$$ 
   is  fully faithful.
   All that is left to show is that it is
   essentially surjective. Since $\grPGL(n_0,n_1,\cdots,n_r)$ and
   $\AUT\PP(n_0,n_1,\cdots,n_r)$ are both stacks, it is enough to
   prove this for $U=\Spec R$, where $R$ is a local ring.
   In this case, we know  by Proposition \ref{P:linebundles} that
   $\Pic\PP(n_0,n_1,\cdots,n_r)\cong\bbZ$.
   We can now proceed  exactly as in (\oldcite{B-N}, Theorem
   8.1).

    The isomorphisms stated at the end of the theorem follow
    from Lemma \ref{L:isom} and Lemma
    \ref{L:easy4}.
\end{proof}

\section{Structure of $PGL(n_0,n_1,...,n_r)$}{\label{S:Structure}}

In this section we give detailed information about the structure of
the group $G_{n_0,n_1,\cdots,n_r}$. We show that, as a group scheme
over an arbitrary base, it splits as a semi-direct product of a
reductive group scheme and a unipotent group scheme. The reductive
part is a product of a copies of the general linear groups. The
unipotent part is a successive semi-direct product of vector groups;
see Theorem \ref{T:decomposition}.

  Throughout this section, the action of $\Gm$ on $\bbU=\Ao{r+1}$
  means the weight $(n_0,n_1,\cdots,n_r)$ action. To shorten the
  notation,
  we denote the group $G_{n_0,n_1,\cdots,n_r}$
  by $G$. The rank $m$ general linear group scheme over $\Spec R$
  is denoted by $\GL(m,R)$. When $R=\bbZ$, this is abbreviated to
  $\GL(m)$. We always assume $r\geq 1$. The corresponding
  projectivized group scheme is denoted by $\PGL(m)$;
  this notation does not conflict with the notation $\PGL(n_0,n_1,\cdots,n_r)$
  for a  weighted projective general linear 2-group as in the
  latter case we have at least two variables.

 We begin with a simple lemma.

\begin{lem}{\label{L:extend}}
    Let $R$ be an arbitrary ring, and let $f$ be a global section
    of the structure sheaf of  $\bbU_R=\bbA^{r+1}_R-\{0\}$,  $r \geq 1$.
    Then $f$ extends uniquely to a global section of $\bbA^{r+1}_R$.
\end{lem}

\begin{proof}
    Let $U_i=\Spec R[x_0,\cdots,x_r,x_i^{-1}]$ and  consider the
    covering $\bbU_R=\cup_{i=1}^{n} U_i$. We show that the
    restriction $f_i:=f|_{U_i}$ is a polynomial for every $i$.
    To see this, observe that, except possibly for $x_i$, all
    variables occur with positive powers in $f_i$. To show that $x_i$ also
    occurs with a positive power, pick some $j\neq i$ and use the
    fact that $x_i$ occurs with a positive power in
    $f_j|_{U_i\cap U_j}=f_i|_{U_i\cap U_j}$.

    Therefore,  for every $i$, $f_i$ actually lies in $R[x_0,\cdots,x_r,x_i^{-1}]$.
    Since  $f_j|_{U_i}=f_i|_{U_j}$, it is obvious that all $f_i$
    are actually the same and provide the desired extension of $f$
    to $\bbU_R$.
 \end{proof}

  From now on, we will use a slightly
  different notation with indices. Namely, we assume that the weights are
  $m_1<m_2<\cdots<m_t$, with each $m_i$ appearing exactly
  $r_i \geq 1$ times in the weight sequence (so in the previous notation
  we would have $r+1=r_1+\cdots+r_t$). We denote the corresponding
  projective general linear 2-group by $\PGL(m_1:r_1,m_2:r_2,\cdots,m_t:r_t)$.
  We use the coordinates $x^i_j$, $1\leq i\leq t$, $1\leq j\leq
  r_i$,   for $\bbA^{r+1}$. We think of $x^i_j$ as a variable of degree $m_i$.
 We will usually  abbreviate the sequence  $x^i_1,\cdots,x^i_{r_i}$ to $\mathbf{x}^i$. 
 Similarly,  a sequence  $F^i_1,\cdots,F^i_{r_i}$ of  polynomials is abbreviated to
 $\mathbf{F}^i$.

\vspace{0.1in}

  Let $R$ be a ring. The following proposition tells us how a
  $\Gmr{R}$-equivariant automorphisms of   $\bbU_R$ looks like.

  \begin{prop}{\label{P:endomorphisms}}
   Let $F \: \bbU_R \to \bbU_R$  be a $\Gm$-equivariant map. Then $F$ is of the form
  $(\mathbf{F}^i)_{1\leq i\leq t}$, where for every $i$,
   each component $F^i_j \in R[x^i_j; 1\leq i\leq t,1\leq j\leq r_i]$ of
    $\mathbf{F}^i$ is a weighted   homogeneous   polynomial of weight $m_i$.
  \end{prop}

 \begin{proof}
   The fact that components of $F$ are polynomial  follows from
   Lemma \ref{L:extend}.   The statement about   homogeneity of $F^i_j$ is a simple exercise
   in polynomial algebra and is left to the reader.
 \end{proof}

 In the above proposition, each $F^i_j$ can be written in the form
 $F^i_j=L^i_j+P^i_j$, where $L^i_j$ is linear in the variables
 $x^i_1,\cdots, x^i_{r_i}$, and $P^i_j$ is a homogeneous polynomial
 of degree $m_i$ in variables $x^a_b$ with $a < i$. Let
 $L_F:=(\mathbf{L}^i)_{1\leq i\leq t}$ be the linear part of $F$.
 It is again a $\Gm$-equivariant endomorphism of $\bbU$.

 \begin{prop}{\label{P:linear}}
   Let $F$ be as in the Proposition \ref{P:endomorphisms}.
   The assignment $F  \mapsto L_F$ respects composition of
   endomorphisms.  In particular, if $F$ is an automorphism, then so is $L_F$.
 \end{prop}

 \begin{proof}
    This follows from direct calculation, or, alternatively, by using the fact that
    $L_F$ is simply the derivative of $F$ at the origin.
 \end{proof}

 \begin{cor}{\label{C:linear}}
    There is a natural split homomorphism
       $$\phi \: G \to \GL(r_1)\times\GL(r_{2})\times\cdots\times\GL(r_t).$$
 \end{cor}

 Next we give some information about the structure of the
 kernel $U$  of $\phi$. It consists of endomorphisms
 $F=(F^i_j)_{i,j}$, where $F^i_j$ has the form
   $$F^i_j=x^i_j+P^i_j.$$
 Here, $P^i_j$ is a homogeneous polynomial
 of degree $m_i$ in variables $x^a_b$ with $a < i$.
 Indeed, it is easily seen that, for an arbitrary
 choice of the polynomials $P^i_j$, the resulting endomorphism
 $F$ is automatically invertible. So, to give such an
 $F \in U$  is equivalent to giving an arbitrary collection
 of polynomials $\{P^i_j\}_{1\leq i\leq t,1\leq j\leq r_i}$
 such that each $P^i_j$ is a homogeneous polynomial
 of degree $m_i$ in variables $x^a_b$ with $a < i$.
 So, from now on we switch the notation and
 denote such an element of $U$ by $(P^i_j)_{i,j}$.

 \begin{prop}{\label{P:decomposition}}
  For each $1\leq a\leq t$, let $U_a \subseteq U$ be the set of
  those endomorphisms $F=(P^i_j)_{i,j}$ for which $P^i_j=0$ whenever
  $i\neq a$. Let $K_a$ denote the set of monomials of degree $m_{a}$ in variables
  $x^i_j$, $i < a$, and let $k_a$ be the cardinality of $K_a$.
  (In other words, $k_a$ is the number of solutions of the equation
    $$\sum_{i=1}^{a-1}m_i\sum_{j=1}^{r_i} z_{i,j} = m_a$$
  in non-negative integers $z_{i,j}$.) Then we have the following:

   \begin{itemize}
    \item[($i$)]
      $U_a$ is a subgroup of $U$ and  is
     canonically isomorphic to the vector group scheme
     $\mathbb{A}^{r_a}\otimes\mathbb{A}^{K_a}\cong\mathbb{A}^{r_{a}k_a}$.
    (Note: $U_1$ is trivial.)

    \item[($ii$)] If $a< b$, then $U_a$ normalizes $U_{b}$.

    \item[($iii$)] The groups $U_a$, $1\leq i\leq t$, generate $U$
     and we have $U_a\cap U_{b}=\{1\}$ if $a\neq b$.
  \end{itemize}
 \end{prop}

 \begin{proof}[Proof of ($i$)]

   The action of $(P^i_j)_{i,j}\in U_a$ on $\mathbb{A}^{r+1}$
   is given by
      $$(\mathbf{x}^1,\cdots,\mathbf{x}^a,\cdots,\mathbf{x}^t)
            \longmapsto  (\mathbf{x}^1,\cdots,
                               \mathbf{x}^a+\mathbf{P}^a,\cdots,\mathbf{x}^t).$$
    So, if $\mathbb{A}^{K_a}$ stands for the vector group scheme
    on the basis $K_a$, there is a canonical isomorphism
      $$U_a\cong\bigoplus_{i=1}^{r_a}\mathbb{A}^{K_a}
                             \cong\mathbb{A}^{r_a}\otimes\mathbb{A}^{K_a}.$$

\vspace{0.1in}

 \noindent{\em Proof of} ($ii$). Let $G=(Q^i_j)_{i,j}$
   be an element in $U_a$ and $F=(P^i_j)_{i,j}$ an element in
   $U_{b}$. By ($i$), the inverse of $G$ is
   $G^{-1}=(-Q^i_j)_{i,j}$.
   Let us analyze the effect of the composite $G\circ F\circ G^{-1}$ on
   $\mathbb{A}^{r+1}$:

$$\begin{array}{rcl}
      (\mathbf{x}^1,\cdots,\mathbf{x}^a,\cdots,\mathbf{x}^{b},\cdots,\mathbf{x}^t)
          & \stackrel{G^{-1}}\longmapsto & (\mathbf{x}^1,\cdots,
                               \mathbf{x}^a-\mathbf{Q}^a,\cdots,\mathbf{x}^{b},\cdots,\mathbf{x}^t) \\
          & \stackrel{F}\longmapsto & (\mathbf{x}^1,\cdots,
                  \mathbf{x}^a-\mathbf{Q}^a,\cdots,\mathbf{x}^{b}+\mathbf{R}^{b},\cdots,\mathbf{x}^t) \\
          & \stackrel{G}\longmapsto  & (\mathbf{x}^1,\cdots,
                    \mathbf{x}^a,\cdots,\mathbf{x}^{b}+\mathbf{R}^{b},\cdots,\mathbf{x}^t).
     \end{array}$$
  Here the polynomial $R_k^{b}$, $1\leq k \leq r_{b}$,  is obtained from $P_k^{b}$
  by substituting the variables $x_j^{a}$ with the polynomial $x_j^a-Q_j^{a}$.

\vspace{0.1in}

 \noindent{\em Proof of} ($iii$). Easy.
 \end{proof}

 Part ($ii$)   implies that each $U_a$
 acts by conjugation on each of $U_{a+1}$,
 $U_{a+2}$,$\cdots$, $U_t$.\footnote{All group actions in this section are assumed to be on
the left.}
 To fix the notation, in what follows we let the conjugate of an automorphism $f$
 by an automorphism $g$ to be $g\circ f\circ g^{-1}$.

 \vspace{0.1in}
  \noindent{\bf Notation.} Let $\{U_a\}_{a=1}^t$ be a family of subgroups of a
   group $U$ which satisfies the following properties: 1) Each $U_a$ normalizes
   every $U_{b}$ with $a<b$; \ 2) No two distinct $U_a$
   intersect; \ 3) The $U_a$ generate $U$. In this case, we say that $U$ is a
   successive semi-direct product of the $U_a$, and use the notation
   $U \cong U_t\rtimes \cdots\rtimes U_2\rtimes U_1$.
 \vspace{0.1in}

 The following is an immediate corollary of Proposition \ref{P:decomposition}.

 \begin{cor}{\label{C:decomposition1}}
      There is a natural decomposition of $U$ as a semi-direct product
      $$U \cong U_t\rtimes \cdots\rtimes U_2\rtimes U_1,$$
      where $U_a\cong\mathbb{A}^{r_{a}k_a}$ is the group introduced in
      Proposition \ref{P:decomposition}. (Note that $U_1$ is
      trivial.)
 \end{cor}

 In the next theorem we use the notation  $\mathbb{A}^m$ for
  two things. One that has already appeared is the affine group scheme
  of dimension $m$. When there is a group scheme $G$ involved, we
  also use the notation $\mathbb{A}^m$ for the trivial
  representation of $G$ on $\mathbb{A}^m$.

 \begin{thm}{\label{T:decomposition}}
   There is a natural decomposition of $G$ as a semi-direct product
      $$G \cong U_t\rtimes \cdots\rtimes U_2\rtimes U_1 \rtimes
\big(\GL(r_1)\times\cdots\times\GL(r_t)\big),$$
      where $U_a\cong\mathbb{A}^{r_{a}k_a}$ and $k_a$ is as in
      Proposition \ref{P:decomposition}. (Note that $U_1$ is trivial.)
      Furthermore,
      for every $1\leq a\leq t$,
      the action of $\GL(r_a)$ leaves each $U_{b}$ invariant.
      We also have the following:

  \begin{itemize}

   \item[($i$)] When $a > b$ the induced
         action of $\GL(r_{a})$ on $U_{b}$ is trivial.

   \item[($ii$)] When $a=b$ the induced
         action of $\GL(r_a)$ on $U_a$ is
         naturally isomorphic to the representation $\rho\otimes\mathbb{A}^{K_a}$,
         where $\rho$ is the standard representation of
         $\GL(r_a)$ and $K_a$ is as in
         Proposition \ref{P:decomposition}.
         (Recall that $U_a$ is canonically isomorphic to
         $\mathbb{A}^{r_a}\otimes\mathbb{A}^{K_a}$.)

   \item[($iii$)] When $a < b$  the action of
         $\GL(r_a)$ on $U_{b}$ is naturally isomorphic to the
         representation
         $$\bigoplus_{0\leq l\leq \lfloor\frac{m_{b}}{m_a}\rfloor}
                                    \bbA^{r_{b} d_l} \otimes\hat{\rho}^{\otimes l}.$$
         Here $\hat{\rho}$ stands for the inverse transpose of $\rho$,
         and $d_l$ is the number of monomials of degree $m_{b}$ in variables
         $x^i_j$, $i < b$, $i\neq a$; so $d_l$ also depends on $a$ and $b$.
    (In other words, $d_l$ is the number of      solutions of the equation
    $$\sum_{\substack{i=1 \\ i\neq a}}^{b-1}m_i\sum_{j=1}^{r_i} z_{i,j} = m_{b}-lm_a$$
    in non-negative integers $z_{i,j}$.)
   \end{itemize}

 \end{thm}

 \begin{proof}
  Let $g \in \GL(r_{a})$ and $F \in U_{b}$. As in the proof of
  Proposition \ref{P:decomposition}.$i$, we analyze the
  effect of the composite $g\circ F\circ g^{-1}$ on $\mathbb{A}^{r+1}$.
  The element  $g \in \GL(r_{a})$
   acts on $\mathbb{A}^{r+1}$ as follows: it leaves every component $x_i^j$ invariant if
   $i\neq a$ and on the coordinates $x_1^a,\cdots,x_{r_a}^a$ it acts
   linearly (like the action of an $r_a\times r_a$ matrix  on a column vector).

\vspace{0.1in}

\noindent{\em Proof of} ($i$). The effect of
  $g \in \GL(r_{a})$ only involves the variables
  $x_1^a,\cdots,x_{r_a}^a$ and does not see any other variable,
  whereas the effect of $F \in U_{b}$ only involves the variables
  $x^i_j$, $i \leq b$. Since $b < a$, these two are independent
  of each other. That is, $F$ and $g$ commute.

\vspace{0.1in}

\noindent{\em Proof of} ($ii$).
 Assume $F=(P^i_j)_{i,j}$; so $P^i_j=0$ if $i\neq a$.
 The effect of $g\circ F\circ g^{-1}$ can be described as follows:
  $$\begin{array}{rcl}
      (\mathbf{x}^1,\cdots,\mathbf{x}^a,\cdots,\mathbf{x}^t)
          & \stackrel{g^{-1}}\longmapsto & (\mathbf{x}^1,\cdots,
                               \mathbf{y}^a,\cdots,\mathbf{x}^t) \\
          & \stackrel{F}\longmapsto & (\mathbf{x}^1,\cdots,
                  \mathbf{y}^a+\mathbf{P}^a,\cdots,\mathbf{x}^t) \\
          & \stackrel{g}\longmapsto  & (\mathbf{x}^1,\cdots,
                    \mathbf{x}^a+\mathbf{Q}^a,\cdots,\mathbf{x}^t).
     \end{array}$$

  Here,  $y_j^a$ is the linear combination of $x_1^a,\cdots,x_{r_a}^a$, the
  coefficients being the entries of the $j^{th}$ row of the matrix
  $g^{-1}$. Similarly, $Q_j^a$ is the linear combination of
  $P_1^a,\cdots,P_{r_a}^a$,  coefficients being the entries of the $j^{th}$ 
  row of the matrix  $g$.

\vspace{0.1in}

 \noindent{\em Proof of} ($iii$).
   Assume $F=(P^i_j)_{i,j}$; so $P^i_j=0$ if $i\neq b$. Let
   $\mathbf{y}^a$ be as in ($ii$). The effect of
   $g\circ F\circ g^{-1}$ can be described as follows:
  $$\begin{array}{rcl}
      (\mathbf{x}^1,\cdots,\mathbf{x}^a,\cdots,\mathbf{x}^{b},\cdots,\mathbf{x}^t)
          & \stackrel{g^{-1}}\longmapsto & (\mathbf{x}^1,\cdots,
                               \mathbf{y}^a,\cdots,\mathbf{x}^{b},\cdots,\mathbf{x}^t) \\
          & \stackrel{F}\longmapsto & (\mathbf{x}^1,\cdots,
                  \mathbf{y}^a,\cdots,\mathbf{x}^{b}+\mathbf{R}^{b},\cdots,\mathbf{x}^t) \\
          & \stackrel{g}\longmapsto  & (\mathbf{x}^1,\cdots,
                    \mathbf{x}^a,\cdots,\mathbf{x}^{b}+\mathbf{R}^{b},\cdots,\mathbf{x}^t).
     \end{array}$$
  Here the polynomials $R_k^{b}$, $1\leq k \leq r_{b}$,  are obtained from $P_k^{b}$
  by substituting the variable $x_j^{a}$ with $y_j^{a}$.

  Let $\lambda$ be  the representation of $\GL(r_a)$ on the space
  $V$ of homogenous polynomials of degree $m_{b}$ which acts as follows:
  it takes a  polynomial $P \in V$ and substitutes the
  variables $x_j^{a}$, $1\leq j \leq r_a$, with $y_j^{a}$. From the
  description above, we see that the representation of $\GL(r_a)$
  on $U_{b}$ is a direct sum of $r_{b}$ copies of $\lambda$.
  We will  show that
     $$\lambda \cong \bigoplus_{0\leq l\leq \lfloor\frac{m_{b}}{m_a}\rfloor}
                                    \bbA^{d_l} \otimes\hat{\rho}^{\otimes l}.$$
  To obtain the above decomposition, simply note that a polynomial in $V$
  can be uniquely written in the form
    $$\sum_{0\leq l\leq \lfloor\frac{m_{b}}{m_a}\rfloor} S_lT_l,$$
    where $T_l$ is a homogenous polynomial of degree $lm_a$ in
    variables $x_1^a,\cdots, x_{r_a}^a$, and $S_l$ is a homogenous
    polynomial of degree $m_{b}-lm_a$ in the rest of the
    variables. The action of $\GL(r_a)$ leaves $S_l$ intact and
    acts on  $T_l$ by the $l^{th}$ power of the inverse transpose of
    the standard representation.
 \end{proof}

 The actions of various pieces in the above semi-direct product decomposition,
 though explicit, are tedious  to write down, except for small values of $t$.
 We give some examples in $\S$\ref{S:Examples}.
 
 Let us denote   $U_t\rtimes \cdots\rtimes U_2\rtimes U_1$ by $U$ and define the
 crossed module  
     $$\PGL(n_0,n_1,\cdots,n_r)_{red}:=[\partial \: \Gm \to \GL(r_1)\times\cdots\times\GL(r_t)],$$
where the $k^{th}$ factor of $\partial(\lambda)$ is the size $r_k$  scalar matrix $\lambda^{m_k}$. 
Theorem \ref{T:decomposition} then implies that
      $$\PGL(n_0,n_1,\cdots,n_r)\cong U\rtimes\PGL(n_0,n_1,\cdots,n_r)_{red}.$$
We think of $U$ as the unipotent radical and    $\PGL(n_0,n_1,\cdots,n_r)_{red}$ as the 
reductive part of $\PGL(n_0,n_1,\cdots,n_r)$.

 \begin{rem}
  It is perhaps useful to put the above result in the general
  context of algebraic group theory. Recall that every algebraic
  group  $G$ over a {\em field} fits in a short exact sequence
    $$1 \to U \to G \to G_{red} \to 1,$$
  where $U$ is the unipotent radical of $G$ and $G_{red}$ is
  reductive. The sequence is not split in general. In our case,
  the group scheme $G_{n_0,n_1,\cdots,n_r}$ admits such a short
  sequence over an arbitrary base and, furthermore, the sequence is
  split.

  The general theory of unipotent groups tells us that any unipotent
  group over a {\em perfect field} admits a filtration whose graded
  pieces are vector groups. This filtration splits, but only in the category
  of schemes (i.e., the splitting maps may not be group
  homomorphisms). In our case, however, the group scheme $U$ admits
  such a filtration over an arbitrary base. Furthermore, the
  filtration is split group theoretically.
 \end{rem}

\section{Some examples}{\label{S:Examples}}

In this section we look at some explicit examples of weighted projective general linear 2-groups.

 \begin{ex}{\label{E:1}} {\em Weight sequence} $m < n, m\nmid n$.
   In this case we have $t=2$, and $r_1=r_2=1$ and $k_1=0$. So
   $G \cong \Gm\times \Gm$.
 \end{ex}

 \begin{ex}{\label{E:2}} {\em Weight sequence} $m < n, m\mid n$.
   In this case we have $t=2$, $r_1=r_2=1$, and $k_1=1$.
   So we have
     $$G \cong \bbA \rtimes(\Gm\times\Gm).$$
    The action of an element $(\lambda_1,\lambda_2)
    \in \Gm\times\Gm$
    on an element $a \in \bbA$ is given by
       $$(\lambda_1,\lambda_2)\cdot a=\lambda_2\lambda_1^{-\frac{n}{m}}a.$$
    More explicitly, an element in $G$ is map of the form
       $$(x,y) \mapsto (\lambda_1x, \lambda_2y+ax^{\frac{n}{m}}).$$
    Note the similarity with the group of $2\times 2$ lower-triangular
    matrices.
 \end{ex}

 \begin{ex}{\label{E:3}} {\em Weight sequence} $n=m$.
    We obviously have $G\cong \GL(2)$.
\end{ex}

 \begin{ex}{\label{E:4}} {\em Weight sequence} $1,2,3$.
   First we determine $U$. A typical element in $U$ is
   of the form
      $$(x,y,z) \mapsto (x,y+ax^2,z+bx^3+cxy).$$
  We have $U_2=\bbA$ and $U_3=\bbA^2$.
  The action of an element $a \in U_2$ on an
  element $(b,c) \in  U_3$ is given by $(b-ac,c)$.
  That is, $a$ acts on $U_3=\bbA^2$ by the matrix
     $$\left( \begin{array}{cc} 1 & -a \\
                             0 & 1 \end{array}\right) $$
  So, $U\cong \bbA^{\oplus 2}\rtimes \bbA$.
  Finally, we have
     $$G\cong U \rtimes (\Gm)^3=\bbA^{\oplus 2}\rtimes \bbA \rtimes (\Gm)^3,$$
  where the action of an element   $(\lambda_1,\lambda_2,\lambda_3)
  \in (\Gm)^3$
  on an element $(a,b,c) \in U$ is given by
  $(\lambda_1^{-2}\lambda_2a,\lambda_1^{-3}\lambda_3b,
  \lambda_1^{-2}\lambda_2^{-1}\lambda_3c)$.
 \end{ex}

 \begin{ex}{\label{E:5}} {\em Weight sequence} $1,2,4$.
    An element in $U$ has the general form
      $$(x,y,z) \mapsto (x,y+ax^2,z+bx^4+cx^2y+dy^2).$$
  We have $U_2=\bbA$ and $U_3=\bbA^3$.
  The action of an element $a \in U_2$ on an
  element $(b,c,d) \in  U_3$ is given by the matrix
    $$\left( \begin{array}{ccc} 1 & -a  & a^2 \\
                                0 &  1  & -2a \\
                                0 &  0  &  1  \end{array}\right) $$
  So, $U\cong \bbA^{\oplus 3}\rtimes\bbA$.

  Finally, we have
     $$G\cong U\rtimes (\Gm)^3=\bbA^{\oplus 3}\rtimes\bbA  \rtimes (\Gm)^3,$$
  where the action of an element   $(\lambda_1,\lambda_2,\lambda_3)
            \in (\Gm)^3$
  on an element $(a,b,c,d) \in U$ is given by
      $$(\lambda_1^{-2}\lambda_2a,\lambda_1^{-4}\lambda_3b,
          \lambda_1^{-2}\lambda_2^{-1}\lambda_3c,\lambda_2^{-2}\lambda_3d).$$
 \end{ex}

\vspace{0.1in}

 Next we look at
$\PGL(m_1:r_1,m_2:r_2,\cdots,m_t:r_t)$. Recall that, as a
crossed module, this is given by $[\partial \: \Gm \to G]$, where
$\partial$ is the obvious map coming from the action of $\Gm$ on
$\mathbb{A}^{r+1}$, and the action of $G$ on $\Gm$ is the trivial
one.

Observe that the map $\partial$ factors though the component
$\GL(r_1)\times\cdots\times\GL(r_t)$ of $G$. So, let us define $L$
to be the cokernel of the following map:
   $$\xymatrix@=12pt@M=10pt{\Gm
      \ar[rrrrr]^(0.39){\small{
            (\overbrace{\lambda^{m_1},\ldots,\lambda^{m_1}}^{r_1},
                    \cdots\cdot\cdot,
                \overbrace{\lambda^{m_t},\ldots,\lambda^{m_t}}^{r_t})} } &&&&&
                                       \GL(r_1)\times\dots\times\GL(r_t).   }$$
From Theorem \ref{T:decomposition} we immediately obtain the
following.

\begin{prop}{\label{P:PGL}}
   Let $L$ be the group defined in the previous paragraph, and let
   $k=\gcd(m_1,\cdots,m_t)$.
   We have natural isomorphisms of group schemes
     $$\begin{array}{rcl}
          \pi_0\PGL(m_1:r_1,m_2:r_2,\cdots,m_t:r_t)  &\cong&
                    U_t\rtimes \cdots\rtimes U_2\rtimes U_1  \rtimes L, \\
          \pi_1\PGL(m_1:r_1,m_2:r_2,\cdots,m_t:r_t)  &\cong&    \mu_k.
        \end{array}$$
\end{prop}

 Our final result is that, if all weights are distinct (that is, $r_i=1$),
 then the corresponding projective general linear 2-group is split.

 \begin{prop}{\label{P:split}}
   Let $\{m_1,\cdots,m_t\}$ be distinct positive integers, and consider
   the projective general linear 2-group $\PGL(m_1,m_2,\cdots,m_t)$.
   Then, the projection map $G \to \pi_0\PGL(m_1,m_2,\cdots,m_t)$ splits.
   In particular, $\PGL(m_1,m_2,\cdots,m_t)$ is split. That is, it is
   completely classified by its homotopy group schemes:
      $$\begin{array}{rcl}
          \pi_0\PGL(m_1,\cdots,m_t)  &\cong&
              U_t\rtimes \cdots\rtimes U_2\rtimes U_1  \rtimes(\Gm)^{t-1}, \\
          \pi_1\PGL(m_1,\cdots,m_t)  &\cong&    \mu_k.
        \end{array}$$
 \end{prop}

 \begin{proof}
   By Theorem \ref{T:decomposition} and  Proposition \ref{P:PGL} we
   know that
   $G\cong U_t\rtimes \cdots\rtimes U_2\rtimes U_1  \rtimes(\Gm)^{t}$
   and $\pi_0\PGL(m_1,m_2,\cdots,m_t)  \cong
                    U_t\rtimes \cdots\rtimes U_2\rtimes U_1  \rtimes L$, where
    $L$ is the cokernel of the map
     $$\xymatrix@=21pt@M=6pt{\alpha \: \Gm
      \ar[rr]^(0.45){(\lambda^{m_1},\cdots,\lambda^{m_t})}  &&
                                                (\Gm)^t.   }$$
    So it is enough to show that the image of $\mu$ is a direct factor.
    Note that if we divide all the $m_i$ by their greatest common divisor,
    the image of $\alpha$ does not change. So, we may assume
    $\gcd(m_1,\cdots,m_t)=1$. Let $M$ be a $t\times t$ integer matrix
    whose determinant is $1$ and whose first column is $(m_1,\cdots,m_t)$.
    The matrix $M$ gives rise to an isomorphism
    $\mu \: (\Gm)^t \to (\Gm)^t$ whose restriction to the
    subgroup $\Gm\times\{1\}^{t-1}\cong\Gm$ is naturally identified with
    $\alpha$. The subgroup
    $\mu(\{1\}\times (\Gm)^{t-1}) \subset (\Gm)^t$ is
    the desired complement of the image of $\alpha$.
 \end{proof}

 \begin{cor}{\label{C:split2}}
   Let $m$, $n$ be distinct positive integers, and let $k=\gcd(m,n)$.
   Then $\PGL(m,n)$ is a split 2-group. That is, it is classified by its
   homotopy groups:
       $$\begin{array}{rcl}
          \pi_0\PGL(m,n)  &\cong&
                    \left\{  \begin{array}{l}
                        \Gm,  \ \text{if}\ m < n, m\nmid n \\
            \bbA\rtimes\Gm,  \ \text{if}\ m < n, m\mid n
                                                        \end{array}  \right.  \\
           \pi_1\PGL(m,n)  &\cong&    \mu_k.
        \end{array}$$
    (In the case $m\mid n$, the action of $\Gm$ on $\bbA$
    in the cross product $\bbA\rtimes\Gm$ is simply the
    multiplication action.)
 \end{cor}

 \begin{proof}
    Everything is clear, except perhaps a clarification is in order
    regarding the parenthesized statement. Observe that the $\Gm$
    appearing in the  cross product  $\bbA\rtimes\Gm$ is
    indeed the cokernel of the map
        $$\xymatrix@=10pt@M=6pt{\alpha \: \Gm
      \ar[rr]^(0.46){(\lambda^{m},\lambda^{n})}  &&
                                                (\Gm)^2,   }$$
    which is naturally identified  with the subgroup
    $\{1\}\times\Gm \subset (\Gm)^2$. Therefore, by the
    formula of Example \ref{E:2}, the action of an element
    $\lambda  \in \Gm$
    on an element $a \in \bbA$ is given by $\lambda a$.
 \end{proof}

   Finally, for the sake of completeness, we include the following.

  \begin{prop}{\label{P:nonsplit}}
     The 2-group $\PGL(k,k,\cdots,k)$, $k$ appearing $t$ times,
     is given by the following crossed module:
         $$\xymatrix@=10pt@M=6pt{[\Gm
      \ar[rr]^(0.46){(\lambda^{k},\cdots,\lambda^{k})}  &&
                                                \GL(t)].   }$$
    We have $\pi_0\PGL(k,\cdots,k)\cong \PGL(t)$ and  $\pi_1\PGL(k,\cdots,k)\cong \mu_k$.
    In particular, $\PGL(1,1,\cdots,1)$, $1$ appearing $t$ times, is
    equivalent to the group scheme $\PGL(t)$.
  \end{prop}

\section{2-group actions on weighted projective stacks}{\label{S:Application}}

In this section we combine the method developed in  $\S$\ref{S:2gp}--$\S$\ref{S:Action} 
with the results about the structure of $\AUT\PP(n_0,n_1,\cdots,n_r)$  
to study 2-group actions on a weighted projective stack 
$\PP(n_0,n_1,\cdots,n_r)$. The goal is to illustrate how
one can classify 2-group actions using butterflies and how to describe the 
corresponding quotient  2-stacks.  Below, all group scheme are assumed to 
be flat and of finite presentation over a fixed base $S$.

Let $\HH$ be a group stack and $[\psi \: H_1\to H_0]$ a presentation
for it as a crossed module in schemes. By Theorem \ref{T:weakmorphism}, 
to give an action of $\HH$ on $\PP(n_0,n_1,\cdots,n_r)$ is equivalent to giving a 
butterfly diagram
    $$\xymatrix@C=8pt@R=6pt@M=6pt{ H_1 \ar[rd]^{\kappa} \ar[dd]_{\psi}
                          & & \Gm \ar[ld]_{\iota} \ar[dd]^{\partial} \\
                            & E \ar[ld]^{\sigma} \ar[rd]_(0.35){\rho}  & \\
                                       H_0 & & G_{n_0,n_1,\cdots,n_r}       }$$
 In other
words, to give an action of $\HH$ on $\X$, we need to find a central extension
  $$ 1\to \Gm \to E \to H_0 \to 1$$
of $H_0$ by $\Gm$, together with  

\begin{itemize}
   \item a lift $\kappa \: H_1 \to E$ of  $\psi$ to $E$ such that
     $\kappa(\beta^{\sigma(x)})=x^{-1}\kappa(\beta) x$, for every 
     $\beta \in H_1$ and $x \in E$;
   
   \item an extension of the weighted action
   of $\Gm$ to a linear action of $E$ on $\bbA^{r+1}$ which is trivial on
   the image of $\kappa$.
\end{itemize} 

The following result is more or less immediate from the above description 
of a group action.

\begin{thm}{\label{T:lift}}
 Let $k$ be a field and $H$ a connected linear algebraic group over $k$, 
 assumed to be reductive  if $\operatorname{char}k > 0$. Let $\X=\PP(n_0,n_1,\cdots,n_r)$
 be a wighted projective stack over $k$. Suppose that  $\operatorname{Pic}(H)=0$.
 Then, every action of $H$ on $\X$ lifts to an action of $H$ on $\bbA^{r+1}$ via a 
 homomorphism $H \to G_{n_0,n_1,\cdots,n_r}$.
\end{thm}

\begin{proof}
   An action of $H$ on $\X$ lifts to $\bbA^{r+1}$ if and only if the corresponding 
   butterfly is equivalent to a strict   one. By (\cite{ButterflyI}, Proposition 4.5.3) 
   this is equivalent to the central extension
          $$ 1\to \Gm \to E \to H \to 1$$
   being split. By (\cite{Colliot-Thelene}, Corollary 5.7) such central extensions are classified
   by  $\operatorname{Pic}(H)$, which in our case is assumed to be trivial.  Any choice of a
   splitting amounts to a lift of the action of $H$ to  $\bbA^{r+1}$.
\end{proof}

This result is essentially saying that  the obstruction to lifting the  
$H$-action from  $\PP(n_0,n_1,\cdots,n_r)$ to $\bbA^{r+1}$ lies in 
$\operatorname{Pic}(H)$.

\subsection{Description of the quotient 2-stack}{\label{SS:2stack}}

Given an action of a group stack $\HH$ on the weighted 
projective stack $\PP(n_0,n_1,\cdots,n_r)$, we can use the associated butterfly to
get information about the quotient 2-stack $[\PP(n_0,n_1,\cdots,n_r)/\HH]$. First, notice that
$[\kappa \: H_1 \to E]$ is also a crossed module in schemes. Denote the associated
group stack by $\HH'$. It acts on  $\PP(n_0,n_1,\cdots,n_r)$ via $\rho$. We have
     $$[\PP(n_0,n_1,\cdots,n_r)/\HH]\cong[(\bbA^{r+1}-\{0\})/\HH'].$$
Note that the right hand side is the quotient stack of the action of a group stack on an honest scheme,
namely, $\bbA^{r+1}-\{0\}$. It is easy  to describe its 2-stack structure by looking at
cohomologies of the NW-SE sequence
  $$H_1 \llra{\kappa} E \llra{\rho}  G_{n_0,n_1,\cdots,n_r} $$
of the butterfly. Set
      $$[\PP(n_0,n_1,\cdots,n_r)/\HH]_1:=[(\bbA^{r+1}-\{0\})/\coker\kappa],$$
 and      
    $$[\PP(n_0,n_1,\cdots,n_r)/\HH]_0:=[(\bbA^{r+1}-\{0\})/\im\rho].$$
 (Here, $\coker\kappa$ and $\im\rho$ are the sheaf theoretic cokernel and image of
 the corresponding maps which we assume are representable.)  
Then, $[\PP(n_0,n_1,\cdots,n_r)/\HH]_1$ is the best  approximation of 
$[\PP(n_0,n_1,\cdots,n_r)/\HH]$ by a 1-stack, in the sense that it is obtained by killing off
the 2-automorphisms of the 2-stack $[\PP(n_0,n_1,\cdots,n_r)/\HH]$. 
More precisely,   there is a natural map 
  $$[\PP(n_0,n_1,\cdots,n_r)/\HH] \to [\PP(n_0,n_1,\cdots,n_r)/\HH]_1$$
 making the former a 2-gerbe over the latter for the 2-group $[\ker\kappa \to 1]$.
 Similarly,  $[\PP(n_0,n_1,\cdots,n_r)/\HH]_0$ is an orbifold (i.e., a Deligne-Mumford
 stack which is generically a scheme) obtained by quotienting out the generic 1-automorphisms.
 More precisely, there is a natural map 
  $$[\PP(n_0,n_1,\cdots,n_r)/\HH]_1 \to [\PP(n_0,n_1,\cdots,n_r)/\HH]_0$$
 making the former a gerbe over the latter for the group $\ker\rho/\im\kappa$ (namely, the middle
 cohomology of the NW-SE sequence).

It follows that the quotient 2-stack $[\PP(n_0,n_1,\cdots,n_r)/\HH]$ is equivalent to a 1-stack if and only if
$\kappa$ is injective; it is an orbifold if and only if the NW-SE sequence is left exact.

\begin{ex}{\label{E:case1}}
  Suppose that $\HH$ is an honest  group scheme and denote it by $H$. Then, to give an action
  of $H$ on $\PP(n_0,n_1,\cdots,n_r)$ is equivalent to giving  a central extension
      $$ 1\to \Gm \to E \to H \to 1$$
   of $H$ by $\Gm$, together with a linear action of $E$ on $\bbA^{r+1}$ extending the weighted action
   of $\Gm$. We have
          $$[\PP(n_0,n_1,\cdots,n_r)/H]\cong[(\bbA^{r+1}-\{0\})/E],$$
  which is an honest 1-stack.
\end{ex}

\begin{ex}{\label{E:case2}}
 Suppose that $\HH$ is the group stack associated to $[A\to 1]$, where $A$ is an abelian group scheme.
  We rename $\HH$ to $A[1]$. Then, to give an action of $A[1]$ on 
   $\PP(n_0,n_1,\cdots,n_r)$ is equivalent to giving a character $\kappa \: A \to \mu_d\subset\Gm$,
   where $d$ is the greatest common divisor of $(n_0,n_1,\cdots,n_r)$. The quotient 2-stack
   $[\PP(n_0,n_1,\cdots,n_r)/A[1]]$ is a 1-stack if and only if $\kappa$ is injective. Assume this to be the
   case and identify $A$ with the corresponding subgroup of $\mu_d$. Then,
   roughly speaking, the quotient stack $[\PP(n_0,n_1,\cdots,n_r)/A[1]]$ is obtained by
   killing  the $A$  in  $\mu_d \subseteq I_x$ at every inertia group $I_x$ 
   of $\PP(n_0,n_1,\cdots,n_r)$. (Note that the generic inertia group
   of $\PP(n_0,n_1,\cdots,n_r)$ is $\mu_d$.)
   For example, if the base is an algebraically closed field of 
   characteristic prime to $d$, then 
        $$[\PP(n_0,n_1,\cdots,n_r)/A[1]]\cong \PP(\frac{n_0}{a},\frac{n_1}{a},\cdots,\frac{n_r}{a}),$$
   where $a$ is the order of $A$.
\end{ex}

\providecommand{\bysame}{\leavevmode\hbox
 to3em{\hrulefill}\thinspace}
\providecommand{\MR}{\relax\ifhmode\unskip\space\fi MR }
\providecommand{\MRhref}[2]{%
  \href{http://www.ams.org/mathscinet-getitem?mr=#1}{#2}
} \providecommand{\href}[2]{#2}


\begin{thebibliography}{1}

\bibitem[AlNo1]{ButterflyI} E.~Aldrovandi, B.~Noohi,
  \emph{Butterflies I: morphisms of 2-group stacks}, Adv. Math. 
  \textbf{221} (2009), no. 3, 687--773.

\bibitem[BeNo]{B-N}  K.~Behrend, B.~Noohi, \emph{Uniformization of
{D}eligne-{M}umford analytic
  curves}, J. reine angew. Math. \textbf{599} (2006), 111--153.

\bibitem[Bre]{Breen} L.~Breen, \emph{On the classification of
$2$-gerbes and $2$-stacks}, Ast\'erisque No. 225, 1994.

\bibitem[Bro]{Brochard} S.~Brochard,  \emph{Foncteur de Picard d'un champ alg\'ebrique}, 
Math. Ann. \textbf{343} (2009), no. 3, 541--602.

\bibitem[C-T]{Colliot-Thelene} J-L.~Colliot-Th\'el\`ene, \emph{R\'esolutions flasques des groupes lin\'eaires connexes}, 
 J. reine angew. Math. \textbf{618} (2008), 77--133.

\bibitem[Ha]{Hart}
R.~Hartshorne, \emph{Algebraic geometry}, Graduate Texts in
Mathematics, No.
  52, Springer-Verlag, New York-Heidelberg, 1977.

\bibitem[Hi]{Hirsch} P.~S. Hirschhorn, \emph{Model categories and
their localizations}, Mathematical Surveys and Monographs, 99, American Mathematical Society, Providence, RI,
  2003.

\bibitem[Ho]{Hollander} Sh. Hollander, \emph{A homotopy theory for
stacks}, math.AT/0110247.

\bibitem[No1]{Notes} B.~Noohi, \emph{Notes on 2-groupoids, 2-groups and
crossed modules}, Homotopy, Homology, and Applications, \textbf{9}
(2007), no. 1, 75--106.

\bibitem[No2]{Cohomology} \bysame, \emph{Group cohomology with coefficients in a crossed module}, 
J. Inst. Math. Jussieu,  \textbf{10} (2011), no. 2, 359--404. 

\bibitem[No3]{Maps} \bysame, \emph{On weak maps between 2-groups},
preprint, arXiv:math/0506313v3 [math.CT].

\bibitem[No4]{Picard} \bysame, \emph{Picard stack of a weighted
projective stack}, preprint available at
\texttt{http://www.maths.qmul.ac.uk/$\sim$noohi/research.html}.

\bibitem[Ro]{Romagny}  M.~Romagny, \emph{Group actions on stacks and applications},
Michigan Math. J. \textbf{53}, Issue 1 (2005), 209--236.

\end{thebibliography}
\end{document}